# Evaluation of series with Hurwitz and Lerch zeta function coefficients by using Hankel contour integrals.

Khristo N. Boyadzhiev

**Abstract**. We introduce a new technique for evaluation of series with zeta coefficients and also for evaluation of certain integrals involving the logGamma function. This technique is based on Hankel integral representations of the Hurwitz zeta, the Lerch Transcendent, the Digamma and logGamma functions.

## 1. Introduction.

The Hurwitz zeta function $\zeta(s,a)$ is defined for all $\operatorname{Re} s > 1$, $\operatorname{Re} a > 0$ by

$$\zeta(s,a) = \sum_{n=0}^{\infty} \frac{1}{(n+a)^s}, \qquad (1.1)$$

and has the integral representation:

$$\zeta(s,a) = \frac{1}{\Gamma(s)} \int_0^{+\infty} \frac{x^{s-1}}{1-e^{-x}} e^{-ax} dx. \qquad (1.2)$$

When $a = 1$, it turns into Riemann's zeta function, $\zeta(s,1) = \zeta(s)$.

In this note we present a new method for evaluating the series

$$S(t,a,p) = \sum_{n=1}^{\infty} \frac{\zeta(n+1,a)}{n+p} t^{n+p}, \qquad (1.3)$$

and

$$T(t,a,p) = \sum_{n=1}^{\infty} \frac{\zeta(n+1,a)}{(n+1)(n+2)\ldots(n+p)} t^{n+p} \qquad (1.4)$$

($p = 0, 1, 2, \ldots$) in a closed form. The two series have received a considerable attention since Srivastava [17], [18] initiated their systematic study in 1988. Many interesting results were obtained consequently by Srivastava and Choi (for instance, [6]) and were collected in their recent book [19]. Fundamental contributions to this theory and independent evaluations belong also to Adamchik [1] and Kanemitsu et al [13], [15], [16], Hashimoto et al [12]. For some recent developments see [14].

The technique presented here is very straightforward and applies also to series with the Lerch Transcendent [8]:

$$\Phi(\lambda, s, a) = \sum_{n=0}^{\infty} \frac{\lambda^n}{(n+a)^s}, \tag{1.5}$$

in the coefficients. For example, we evaluate here in a closed form the series

$$\sum_{n=0}^{\infty} \frac{\Phi(\lambda, n+1, a)}{n+p} t^{n+p} \tag{1.6}$$

The evaluation of (1.3) and (1.4) requires zeta values $\zeta(m, a)$ for positive and negative integers $m$. We use a representation of $\zeta(s, a)$ in terms of a Hankel integral, which makes it possible to represent the values $\zeta(m, a)$ for positive and negative integers by the same type of integral. The series (1.3) and (1.4) are evaluated here in terms of the functions

$$g(n, a) = \zeta'(-n, a) + \psi(n+1)\zeta(-n, a), \; n = 0, 1, \ldots, \tag{1,7}$$

and similar expressions are used for (1.6). Here $\psi(s) = \Gamma'(s)/\Gamma(s)$ is the Digamma function and

$$\psi(n+1) = -\gamma + 1 + \frac{1}{2} + \ldots + \frac{1}{n} \tag{1.8}$$

($\gamma$ is Euler's constant). The values $\zeta(-n, a)$ are also well-known and can be written in terms of Bernoulli's polynomials [3], [19]:

$$\zeta(-n, a) = \frac{-B_{n+1}(a)}{n+1}. \tag{1.9}$$

The values $\zeta'(-n, a)$, however, represent a challenge and give rise to new constructions [1], [7],

[10], [14], [19].

In section 2 we introduce Hankel contours and obtain integral representations for $\zeta(m,a)$ and $g(n,a)$. Hankel integral representations for the Digamma function $\psi(s)$ and the logGamma function $\log\Gamma(s)$ are given in section 3. The series (1.3) and (1.4) are evaluated in section 4, while sections 5 and 6 deal with the series (1.6).

Our method makes it possible to evaluate integrals of the form

$$\int_0^t s^m f(a+s)ds, \ Re(a) > 0, \ m = 0,1,2,\ldots,$$

for functions $f(s)$ which have certain Hankel integral representations. Some integrals like

$$\int_0^t s^m \log\Gamma(a+s)ds, \ Re(a) > 0, \ m = 0,1,2,\ldots \qquad (1.10)$$

have already been evaluated by Gosper [11] and Adamchik [1], see also [6], [7], [9], [10], [19]. Using the Hankel integral technique we give an independent evaluation of (1.10) in Section 7. In that section we also obtain a Hankel integral representation for $\log G(s)$, where $G(s)$ is the Barnes $G$-function.

The Appendix at the end of the paper contains a list of some of the Hankel integral representations for easy reference.

.

**2. Hankel integrals**.

For $Re(a) > 0$ consider the integral

$$I(s) = \frac{1}{2\pi i} \int_L \frac{z^{s-1}}{1-e^z} e^{az} dz, \qquad (2.1)$$

where $L$ is the Hankel contour consisting of three parts: $L = L_- \cup L_+ \cup L_\epsilon$, with $L_-$ the "lower side" (i.e. $arg(z) = -\pi$) of the ray $(-\infty, -\epsilon)$, $\epsilon > 0$, traced left to right, and $L_+$ the "upper side" ($arg(z) = \pi$) of this ray traced right to left. Finally, $L_\epsilon = \{z = \epsilon e^{\theta i}: -\pi \leq \theta \leq \pi\}$ is a small circle traced counterclockwise and connecting the two sides of the ray. This contour is used, for example, in [3, p.253] and [20, p.48]. The integral does not depend on $\epsilon$ and setting $\epsilon \to 0$ one can easily see that the contribution from $L_\epsilon$ approaches zero. Evaluating the limit we find:

$$2\pi i\, I(s) = 2\pi i \lim_{\epsilon \to 0} I(s) = \int_{+\infty}^{0} \frac{(xe^{-\pi i})^{s-1} e^{axe^{-\pi i}}}{1 - e^{xe^{-\pi i}}} dx + \int_{0}^{+\infty} \frac{(xe^{\pi i})^{s-1} e^{axe^{\pi i}}}{1 - e^{xe^{\pi i}}} dx$$

$$= -e^{-\pi i s} \int_{0}^{\infty} \frac{x^{s-1} e^{-ax}}{1 - e^{-x}} dx + e^{\pi i s} \int_{0}^{\infty} \frac{x^{s-1} e^{-ax}}{1 - e^{-x}} dx$$

$$= (e^{\pi i s} - e^{-\pi i s}) \int_{0}^{\infty} \frac{x^{s-1} e^{-ax}}{1 - e^{-x}} dx.$$

Therefore,

$$I(s) = \frac{\sin(\pi s)}{\pi} \int_{0}^{\infty} \frac{x^{s-1} e^{-ax}}{1 - e^{-x}} dx. \tag{2.2}$$

In particular, $I(m) = 0$ for every positive integer $m > 1$. From (2.2)

$$\zeta(s,a) = \frac{1}{\Gamma(s)} \frac{\pi}{\sin(\pi s)} I(s), \tag{2.3}$$

and also

$$\zeta(s,a) = \Gamma(1-s) I(s), \tag{2.4}$$

in view of the identity

$$\frac{\pi}{\sin(\pi s)} = \Gamma(s)\Gamma(1-s).$$

Equation (2.4) will be used to represent $\zeta(s,a)$ for negative integer values of $s$, while (2.3) will be used for the positive integer values. First, a note is due: the above equations were obtained for $\mathrm{Re}\, s > 1$, $\mathrm{Re}\, a > 0$; the function $I(s)$, however, is defined and holomorphic for every complex $s$, and therefore equation (2.4) extends and represents the Hurwitz zeta function for every $s \neq 1, 2, \ldots$. Similar representations of the Hurwitz zeta function by a Hankel integral with a contour stretching along the positive semi-axis can be found in [8], [19] and [21].

For $n = 0, 1, 2, \ldots$, (2.4) gives:

$$\zeta(-n, a) = \frac{n!}{2\pi i} \int_L \frac{z^{-(n+1)}}{1 - e^z} e^{az} dz. \tag{2.5}$$

We can not plug $s = n + 1$ in (2.3) or (2.4) directly, so we write (using limit and the rule of L'Hospital in (2.3)):

$$\zeta(n+1, a) = \frac{1}{n!} \{ \lim_{s \to n+1} \frac{\pi I(s)}{\sin(\pi s)} \} = \frac{1}{n!} \{ \frac{I'(n+1)}{\cos(\pi(n+1))} \},$$

$$\zeta(n+1, a) = \frac{(-1)^{n+1}}{2\pi i n!} \int_L \frac{z^n \operatorname{Log} z}{1 - e^z} e^{az} dz. \tag{2.6}$$

The first part of the following lemma is ready.

**Lemma 1.** Equation (2.5) holds for $n = 0, 1, 2, \ldots$ and equation (2.6) holds for $n = 1, 2, \ldots$. Also, for all $s \neq 1, 2, \ldots$

$$\zeta'(s, a) + \psi(1 - s) \zeta(s, a) = \frac{\Gamma(1 - s)}{2\pi i} \int_L \frac{z^{s-1} e^{az} \operatorname{Log} z}{1 - e^z} dz, \tag{2.7}$$

In particular, for $n = 0, 1, 2, \ldots$:

$$g(n, a) \equiv \zeta'(-n, a) + \psi(n+1) \zeta(-n, a) = \frac{n!}{2\pi i} \int_L \frac{z^{-n-1} e^{az} \operatorname{Log} z}{1 - e^z} dz, \tag{2.8}$$

where $\psi(s)$ is the Digamma function (see below section 3).

*Proof.* Equation (2.7) results from (2.4) by differentiation for $s$:

$$\zeta'(s, a) = -\Gamma'(1-s) I(s) + \Gamma(1-s) I'(s)$$

$$= \frac{-\Gamma'(1-s)}{\Gamma(1-s)} \Gamma(1-s) I(s) + \Gamma(1-s) I'(s) = -\psi(1-s) \zeta(s, a) + \Gamma(1-s) I'(s),$$

and (2.8) comes from (2.7) with $s = -n$.

**Remark.** It follows immediately from the representation (2.8) that the functions $g(m,a)$ satisfy the relation:

$$\frac{d}{da}g(m,a) = m\,g(m-1,a), \tag{2.9}$$

or

$$\int_0^t g(m-1, a+s)ds = \frac{1}{m}[g(m,a+t) - g(m,a)]. \tag{2.10}$$

Equation (2.10) can be use used for an alternative proof of Theorem 2 below by consecutive integration.

The functions $g(m,a)$ are similar to the generalized polygamma functions $\psi(-m,q)$ introduced by Espinosa and Moll [9] and the results obtained by them for $\psi(-m,q)$ can easily be written in terms of $g(m,a)$ and vise versa.

### 3. The Digamma function.

At this point we need to provide a Hankel integral representations also for the Gauss Digamma function $\psi(s)$ and the logGamma function $\log\Gamma(s)$. The Digamma function is defined by ([2], [19], [20], [21]):

$$\psi(s) = \frac{d}{ds}\log\Gamma(s) = \frac{\Gamma'(s)}{\Gamma(s)},$$

it has the integral representation (Gauss), [2, p. 26]:

$$\psi(s) = \int_0^\infty \left(\frac{e^{-x}}{x} - \frac{e^{-sx}}{1-e^{-x}}\right)dx, \quad (\operatorname{Re} s > 0), \tag{3.1}$$

and the series representation:

$$\psi(s) = -\gamma + \sum_{k=0}^\infty \left(\frac{1}{k+1} - \frac{1}{k+s}\right), \tag{3.2}$$

where
$$\gamma = -\psi(1) \tag{3.3}$$
is Euler's constant. The values of $\psi$ on the positive integers are given by (1.8).

**Lemma 2.** The Digamma function has the following representations by Hankel integrals $(Re\ s > 0)$:

$$\psi(s) = \frac{1}{2\pi i} \int_L \left(\frac{e^z}{z} + \frac{e^{sz}}{1-e^z}\right) \text{Log } z\, dz, \tag{3.4}$$

also:

$$\psi(s) = \frac{\Gamma(s)}{2\pi i} \int_L z^{-s} e^z \text{Log } z\, dz, \tag{3.5}$$

$$\gamma = -\psi(1) = \frac{-1}{2\pi i} \int_L \frac{e^z \text{Log } z}{z}\, dz, \tag{3.6}$$

$$\psi(s) + \gamma = \frac{1}{2\pi i} \int_L \frac{e^{sz}}{1-e^z} \text{Log } z\, dz. \tag{3.7}$$

*Proof.* Setting $\epsilon \to 0$ in the contour $L$ one computes:

$$\frac{1}{2\pi i} \int_L \left(\frac{e^z}{z} + \frac{e^{sz}}{1-e^z}\right) \text{Log } z\, dz = \frac{-1}{2\pi i} \int_{+\infty}^0 \left(\frac{e^{-x}}{-x} + \frac{e^{-sx}}{1-e^{-x}}\right)(\ln x - \pi i)\, dx$$

$$+ \frac{(-1)}{2\pi i} \int_0^{+\infty} \left(\frac{e^{-x}}{-x} + \frac{e^{-sx}}{1-e^{-x}}\right)(\ln x + \pi i)\, dx + \lim_{\epsilon \to 0} \int_{L_\epsilon} = \psi(s),$$

since the last limit is zero (this is left to the reader). The parts with '$\ln x$' cancel out.

Next, we use the representation of the Gamma function $\Gamma(s)$ by a Hankel contour integral [19, p. 48] (same contour as described above):

$$\frac{1}{\Gamma(s)} = \frac{1}{2\pi i} \int_L z^{-s} e^z\, dz, \tag{3.8}$$

and differentiate both sides for $s$. This leads to (3.5) and (3.6) follows from there when $s = 1$.

Finally, (3.7) is a combination of (3.4) and (3.6). The lemma is proved.

From (3.7), for $Re\ a > 0, Re\ b > 0$ we find

$$\psi(a) - \psi(b) = \frac{1}{2\pi i} \int_L \frac{e^{az} - e^{bz}}{1 - e^z} \text{Log}\, z\, dz, \qquad (3.9)$$

and from (3.2):

$$\psi(a) - \psi(b) = \sum_{k=0}^{\infty} \left( \frac{1}{k+b} - \frac{1}{k+a} \right),$$

which symbolically can be written as

$$\psi(a) - \psi(b) = \zeta(1, b) - \zeta(1, a). \qquad (3.10)$$

In the spirit of Lemma 2 we obtain also a Hankel integral representation for the logGamma function.

**Lemma 3**. The following representation holds:

$$\log \Gamma(a) = \log\sqrt{2\pi} - \gamma\left(a - \frac{1}{2}\right) + \frac{1}{2\pi i} \int_L \frac{z^{-1} e^{az} \text{Log}\, z}{1 - e^z} dz, \qquad (3.11)$$

i.e.

$$\log \Gamma(a) = \log\sqrt{2\pi} - \gamma\left(a - \frac{1}{2}\right) + g(0, a).$$

*Proof.* When $n = 0$ equation (2.8) becomes

$$g(0, a) = \zeta'(0, a) + \psi(1) \zeta(0, a) = \frac{1}{2\pi i} \int_L \frac{z^{-1} e^{az} \text{Log}\, z}{1 - e^z} dz. \qquad (3.12).$$

At the same time, (see, for example, [19], pp. 91-92, or [21], p. 271):

$$\zeta'(0, a) = \log \Gamma(a) - \log\sqrt{2\pi}, \quad \zeta(0, a) = \frac{1}{2} - a. \qquad (3.13)$$

which together with (3.3) lead to (3.11).

**Note.** Differentiating (3.11) for the variable $a$ gives a second proof of (3.7).

### 4. Evaluation of the series with zeta values.

First, we present two simple cases in order to explain the method. Consider the series

$$S = \sum_{n=1}^{\infty} \zeta(n+1, a) t^n , \qquad (4.1)$$

i.e. $S = S(t, a, 0)$. By (2.6) we find

$$S = \frac{1}{2\pi i} \int_L \{\sum_{n=1}^{\infty} \frac{(-1)^{n-1} z^n t^n}{n!}\} \frac{e^{az} \operatorname{Log} z}{1 - e^z} dz,$$

$$S = \frac{1}{2\pi i} \int_L (1 - e^{-tz}) \frac{e^{az} \operatorname{Log} z}{1 - e^z} dz, \qquad (4.2)$$

and therefore, by (3.9):

$$S = \psi(a) - \psi(a - t), \quad (|t| < \operatorname{Re} a). \qquad (4.3)$$

The second example is the series

$$S(t, a, 1) = \sum_{n=1}^{\infty} \frac{\zeta(n+1, a)}{n+1} t^{n+1} , \qquad (4.4)$$

which can be evaluated by integrating (4.3). We shall give a different proof, though, in order to illustrate our method. Instead of (4.3) we integrate (4.2) for $t$ to obtain:

$$S(t, a, 1) = \frac{1}{2\pi i} \int_L (t + z^{-1} e^{-tz} - z^{-1}) \frac{e^{az} \operatorname{Log} z}{1 - e^z} dz \qquad (4.5)$$

(the integration constant is $-z^{-1}$ in order to make the right hand side zero for $t = 0$). Splitting this integral into three parts we write:

$$S(t, a, 1) = \frac{t}{2\pi i} \int_L \frac{e^{az} \operatorname{Log} z}{1 - e^z} dz + \frac{1}{2\pi i} \int_L \frac{z^{-1} e^{(a-t)z} \operatorname{Log} z}{1 - e^z} dz \qquad (4.6)$$

$$-\frac{1}{2\pi i}\int_L \frac{z^{-1}e^{az}\text{Log}\, z}{1-e^z}\,dz,$$

and use (3.7) for the first one and (2.8) for the second and the third:

$$S(t,a,1) = t(\psi(a)+\gamma) + g(0,a-t) - g(0,a).$$

In view of (3.12) and (3.13) this becomes

$$S(t,a,1) = t\psi(a) + \log\Gamma(a-t) - \log\Gamma(a), \quad (|t|<\text{Re}\,a), \qquad (4.7)$$

([21, p. 276]; see also [19, p.159]).

For the general case we need a simple integration formula.

**Lemma 4.** For every $p \in \mathbb{N}$:

$$\int_0^t y^{p-1}(1-e^{-zy})\,dy = \frac{t^p}{p} + e^{-tz}\sum_{k=0}^{p-1} k!\binom{p-1}{k}\frac{t^{p-1-k}}{z^{k+1}} - \frac{(p-1)!}{z^p}. \qquad (4.8)$$

*Proof:* Integration by parts.

We are ready now to evaluate the series (1.3).

**Theorem 1.** For every integer $p \geq 1$, $\text{Re}\,a > 0$, and $|t| < \text{Re}\,a$:

$$S(t,a,p) = \qquad (4.9)$$

$$\frac{t^p}{p}(\psi(a)+\gamma) + \sum_{k=0}^{p-1}\binom{p-1}{k} g(k,a-t)\, t^{p-1-k} - g(p-1,a)$$

where the functions $g(n,a)$ are defined by (2.8).

*Proof.* We multiply both sides of equation (4.2), i.e.

$$\sum_{n=1}^{\infty} \zeta(n+1,a) t^n = \frac{1}{2\pi i}\int_L (1-e^{-tz})\frac{e^{az}\text{Log}\,z}{1-e^z}\,dz$$

by $t^{p-1}$ and integrate for $t$, taking antiderivatives which are zeros at $t=0$. Thus:

$$S(t,a,p) =$$

$$\sum_{n=1}^{\infty} \zeta(n+1,a)\frac{t^{n+p}}{n+p} = \frac{1}{2\pi i}\int_L \{\int_0^t y^{p-1}(1-e^{-yz})\,dy\}\frac{e^{az}\text{Log}\,z}{1-e^z}\,dz. \qquad (4.10)$$

According to (4.8)

$$S(t,a,p) = \frac{t^p}{p}\left\{\frac{1}{2\pi i}\int_L \frac{e^{az}\operatorname{Log} z}{1-e^z}dz\right\} \quad (4.11)$$

$$+ \sum_{k=0}^{p-1}\binom{p-1}{k} t^{p-1-k}\left\{\frac{k!}{2\pi i}\int_L \frac{z^{-k-1}e^{(a-t)z}\operatorname{Log} z}{1-e^z}dz\right\}$$

$$- \frac{(p-1)!}{2\pi i}\int_L \frac{z^{-p}e^{az}\operatorname{Log} z}{1-e^z}dz.$$

The first integral here we replace by $\psi(a)+\gamma$ from (3.7). For the second and the third integrals we use (2.8). This gives (4.9) and the theorem is proved.

In a similar way we can evaluate the series (1.4).

$$T(t,a,p) = \sum_{n=1}^{\infty}\frac{\zeta(n+1,a)}{(n+1)(n+2)\ldots(n+p)}t^{n+p}, \quad (p=1,2,\ldots)$$

**Theorem 2**. For every integer $p \geq 1$, $\operatorname{Re} a > 0$, and $|t| < \operatorname{Re} a$:

$$T(t,a,p) = \frac{t^p}{p!}(\psi(a)+\gamma) \quad (4.13)$$

$$+ \frac{1}{(p-1)!}\sum_{k=0}^{p-1}(-1)^{k+1}\binom{p-1}{k}g(k,a)t^{p-k-1} - \frac{(-1)^p}{p!}g(p,a-t)$$

*Proof.* We integrate equation (4.2) for $t$ consecutively $p$ times to obtain:

$$\sum_{n=1}^{\infty}\frac{\zeta(n+1,a)\,t^{n+p}}{(n+1)(n+2)\ldots(n+p)} = \frac{1}{2\pi i}\int_L Q(t,z,p)\frac{e^{az}\operatorname{Log} z}{1-e^z}dz,$$

where $Q$ and its derivatives for $t$ are zeros at $t=0$, so that

$$Q(t,z,p) = \frac{t^p}{p!} + \sum_{k=0}^{p-1}\frac{(-1)^{k+1}}{(p-k-1)!}z^{-k-1}t^{p-k-1} - (-1)^p z^{-p}e^{-tz}.$$

Therefore,

$$T(t,a,p) = \frac{t^p}{p!} \int_L \frac{e^{az} \operatorname{Log} z}{1-e^z} dz$$

$$+ \sum_{k=0}^{p-1} \frac{(-1)^{k+1} t^{p-k-1}}{(p-k-1)!} \frac{1}{2\pi i} \int_L \frac{z^{-k-1} e^{az} \operatorname{Log} z}{1-e^z} dz - \frac{(-1)^p}{2\pi i} \int_L \frac{z^{-p} e^{(a-t)z} \operatorname{Log} z}{1-e^z} dz.$$

The integrals here can be evaluated by (2.8) and (3.7) and thus one comes to (4.13).

### 5. Series with the Lerch Transcendent

The Lerch Transcendent $\Phi(\lambda, s, a)$, also called Hurwitz-Lerch or Lerch zeta function, is defined by (1.5) for $\operatorname{Re} a > 0$, $\operatorname{Re} s > 0$ and $|\lambda| \leq 1$, $\lambda \neq 1$. When $\lambda = 1$, $\operatorname{Re} s > 1$ it becomes the Hurwitz zeta function $\zeta(s, a)$. We exclude this case which has already been considered. For $a = 1$, the Lerch Transcendent turns into the polylogarithmic function:

$$\lambda \Phi(\lambda, s, 1) = \operatorname{Li}_s(\lambda) = \sum_{m=1}^{\infty} \frac{\lambda^m}{m^s}.$$

Information about $\Phi(\lambda, s, a)$ can be found, for instance, in [8] and [13]. The standard integral representation of this function is similar to (1.2):

$$\Phi(\lambda, s, a) = \frac{1}{\Gamma(s)} \int_0^{+\infty} \frac{x^{s-1}}{1 - \lambda e^{-x}} e^{-ax} dx,$$

and the Hankel integral representations (obtained in the same way as (2.3) and (2.4)) are:

$$\Phi(\lambda, s, a) = \frac{1}{2\pi i \Gamma(s)} \frac{\pi}{\sin(\pi s)} \int_L \frac{z^{s-1} e^{az}}{1 - \lambda e^z} dz, \tag{5.1}$$

$$\Phi(\lambda, s, a) = \frac{\Gamma(1-s)}{2\pi i} \int_L \frac{z^{s-1} e^{az}}{1 - \lambda e^z} dz, \tag{5.2}$$

where $L = L(\epsilon)$ is the same contour described in section 2. Equation (5.2) holds for all $s \in \mathbb{C}$

where both sides are defined. In particular, it is true for $s = 0$. One has for $|\lambda| < 1$:

$$\Phi(\lambda, 0, a) = \sum_{n=0}^{\infty} \lambda^n = \frac{1}{1-\lambda} = \frac{1}{2\pi i} \int_L \frac{z^{-1} e^{az}}{1 - \lambda e^z} dz, \quad (5.3)$$

(note that this integral does not depend on $a$). As before, (5.1) leads to

$$\Phi(\lambda, n+1, a) = \frac{(-1)^{n-1}}{2\pi i \, n!} \int_L \frac{z^n \operatorname{Log} z}{1 - \lambda e^z} e^{az} dz. \quad (5.4)$$

In contrast to (2.6), however, this representation is true also for $n = 0$, i.e.

$$\Phi(\lambda, 1, a) = \frac{-1}{2\pi i} \int_L \frac{e^{az}}{1 - \lambda e^z} \operatorname{Log} z \, dz, \quad (5.5)$$

as both sides are well defined. This is the analog of (3.7). Differentiation of (5.3) gives

$$\Phi_s'(\lambda, -n, a) + \psi(n+1) \Phi(\lambda, -n, a) = \frac{n!}{2\pi i} \int_L \frac{z^{-n-1} e^{az} \operatorname{Log} z}{1 - \lambda e^z} dz, \quad (5.6)$$

which corresponds to (2.8).

**Theorem 3**. For any $|\lambda| \leq 1$, $\lambda \neq 1$, $\operatorname{Re} a > 0$, $|t| < \operatorname{Re} a$ and $p = 1, 2, \ldots$:

$$\sum_{n=0}^{\infty} \Phi(\lambda, n+1, a) \frac{t^{n+p}}{n+p} = \quad (5.7)$$

$$\sum_{k=0}^{p-1} \binom{p-1}{k} [\Phi_s'(\lambda, -k, a-t) + \psi(k+1) \Phi(\lambda, -k, a-t)] t^{p-1-k}$$

$$- [\Phi_s'(\lambda, 1-p, a) + \psi(p) \Phi(\lambda, 1-p, a)].$$

Note that the summation on the left side in (5.7) starts from $n = 0$, while in (4.9) it starts from $n = 1$.

*Proof.* Using (5.4) we find:

$$\sum_{n=1}^{\infty} \Phi(\lambda, n+1, a) t^n = \frac{1}{2\pi i} \int_L (1 - e^{-tz}) \frac{e^{az} \operatorname{Log} z}{1 - \lambda e^z} dz, \quad (5.8)$$

and correspondingly,

$$\sum_{n=1}^{\infty} \Phi(\lambda, n+1, a) \frac{t^{n+p}}{n+p} = \frac{1}{2\pi i} \int_L \{\int_0^t y^{p-1}(1-e^{-yz})\,dy\} \frac{e^{az}\text{Log}\,z}{1-\lambda e^z}\,dz. \qquad (5.9)$$

Proceeding as in the proof of Theorem 1 we come to the analog of (4.11) with the only difference that the denominator in the integrals is $1 - \lambda e^z$ instead of $1 - e^z$. Using (5.5) to substitute $\Phi(\lambda, 1, a)$ for the first integral and (5.6) for the other two, we can write the hand right side of (5.9) as

$$-\Phi(\lambda, 1, a)\frac{t^p}{p} + \sum_{k=0}^{p-1} \binom{p-1}{k} [\Phi'_s(\lambda, -k, a-t) + \psi(k+1)\Phi(\lambda, -k, a-t)]\,t^{p-1-k}$$

$$- [\Phi'_s(\lambda, 1-p, a) + \psi(p)\Phi(\lambda, 1-p, a)].$$

Moving the first term here to the sum on the left side in (5,9) we arrive at (5.7).

When $p = 0$ we have from (5.8)

$$\sum_{n=1}^{\infty} \Phi(\lambda, n+1, a) t^n = \Phi(\lambda, 1, a-t) - \Phi(\lambda, 1, a),$$

or
$$\Phi(\lambda, 1, a-t) = \sum_{n=0}^{\infty} \Phi(\lambda, n+1, a) t^n, \qquad (5.10)$$

which is the Maclaurin series expansion of $\Phi(\lambda, 1, a-t)$ in the variable $t$.

### 6. Evaluating $\Phi(\lambda, -m, a)$ and $\Phi'_s(\lambda, -m, a)$ in terms of more simple functions.

Let $Re(a) > 0$, $|\lambda| < 1$ throughout this section. First we shall evaluate $\Phi(\lambda, -m, a)$ in terms of the geometric polynomials $\omega(\lambda)$ defined by

$$\omega_n(\lambda) = \sum_{k=0}^{n} \{{}^n_k\} k!\,\lambda^k \qquad (6.1)$$

where $\{{}^n_k\}$ are the Stirling numbers of second kind. The geometric polynomials were introduced in [5] and used there for series summation and asymptotic expansions. In particular, the following summation formula holds for every $m = 0, 1, 2, \ldots$

$$(\lambda \frac{d}{d\lambda})^m \frac{1}{1-\lambda} = \sum_{k=0}^{\infty} k^m \lambda^k = \frac{1}{1-\lambda} \omega_m(\frac{\lambda}{1-\lambda}). \tag{6.2}$$

**Proposition 1.** For $m = 0, 1, \ldots$ one has:

$$\Phi(\lambda, -m, a) = \frac{1}{1-\lambda} \sum_{j=0}^{m} \binom{m}{j} a^{m-j} \omega_j(\frac{\lambda}{1-\lambda}) \tag{6.3}$$

*Proof.* When $|\lambda| < 1$ the function

$$\Phi(\lambda, s, a) = \sum_{n=0}^{\infty} \frac{\lambda^n}{(n+a)^s}$$

is defined for all $s \in \mathbb{C}$. In particular, when $s = -m$, $m = 0, 1, 2, \ldots$ one has:

$$\Phi(\lambda, -m, a) = \sum_{n=0}^{\infty} \lambda^n (n+a)^m = \sum_{n=0}^{\infty} \lambda^n \sum_{j=0}^{m} \binom{m}{j} n^j a^{m-j}$$

$$= \sum_{j=0}^{m} \binom{m}{j} a^{m-j} \sum_{n=0}^{\infty} \lambda^n n^j,$$

and according to (6.2) we arrive at (6.3).

Next we consider $\Phi'_s(\lambda, -m, a)$. We shall express these values in terms of the function $l(\lambda, a) = \Phi'_s(\lambda, 0, a)$, i.e.,

$$l(\lambda, a) = -\sum_{n=0}^{\infty} \lambda^n \log(n+a), \tag{6.4}$$

and its derivatives for $\lambda$.

**Proposition 2.** For $m = 0, 1, \ldots$ :

$$\Phi'_s(\lambda, -m, a) = \sum_{q=0}^{m} \binom{m}{q} a^{m-q} (\lambda \frac{d}{d\lambda})^q l(\lambda, a) \tag{6.5}$$

$$= \sum_{q=0}^{m} \binom{m}{q} a^{m-q} \sum_{p=0}^{q} \{{}^q_p\} \lambda^p (\frac{d}{d\lambda})^p l(\lambda,a)$$

*Proof.* One has

$$\Phi'_s(\lambda,s,a) = -\sum_{n=0}^{\infty} \frac{\lambda^n \log(n+a)}{(n+a)^s},$$

and therefore,

$$-\Phi'_s(\lambda,-m,a) = \sum_{n=0}^{\infty} \lambda^n (n+a)^m \log(n+a) = \sum_{q=0}^{m} \binom{m}{q} a^{m-q} \sum_{n=0}^{\infty} \lambda^n n^q \log(n+a). \quad (6.6)$$

Clearly,

$$\sum_{n=0}^{\infty} \lambda^n n^q \log(n+a) = -(\lambda \frac{d}{d\lambda})^q l(\lambda,a). \quad (6.7)$$

At the same time, for any $q$-times differentiable function $f(\lambda)$:

$$(\lambda \frac{d}{d\lambda})^q f(\lambda) = \sum_{p=0}^{q} \{{}^q_p\} \lambda^p f^{(p)}(\lambda), \quad (6.8)$$

(see [5]). Substituting (6.7) in (6.6) we find the representation (6.5) of $\Phi'_s(\lambda,-m,a)$ in terms of $(\frac{d}{d\lambda})^p l(\lambda,a)$.

**Lemma 5.** The function $l(\lambda,a)$ in (6.4) has the integral representation:

$$l(\lambda,a) = \int_0^{\infty} (\frac{e^{-at}}{1-\lambda e^{-t}} - \frac{e^{-t}}{1-\lambda}) \frac{dt}{t}. \quad (6.9)$$

*Proof.* By integrating for $a$ the representation

$$\frac{1}{n+a} = \int_0^{\infty} e^{-(n+a)t} dt,$$

we find that

$$\log(n+a) = \int_0^\infty (e^{-t} - e^{-(n+a)t}) \frac{dt}{t}. \tag{6.10}$$

To obtain (6.9) we multiply (6.10) by $\lambda^n$ and sum from $n = 0$ to $n = +\infty$.

The function $\Phi'_s(\lambda, -m, a)$ has also a representation involving the geometric polynomials in the spirit of (6.3).

**Proposition 3.** For $m = 0, 1, \ldots$:

$$\Phi'_s(\lambda, -m, a) = \sum_{q=0}^{m} \binom{m}{q} a^{m-q} \int_0^\infty [\frac{e^{-at}}{1-\lambda e^{-t}} \omega_q(\frac{\lambda e^{-t}}{1-\lambda e^{-t}}) - \frac{e^{-t}}{1-\lambda} \omega_q(\frac{\lambda}{1-\lambda})] \frac{dt}{t}.$$

This follows from (6.6) where the last sum is replaced by $-(\lambda \frac{d}{d\lambda})^q l(\lambda, a)$, (see (6.7))

and then (6.2) and (6.9) are used. At that, (6.2) is used in the form

$$(\lambda \frac{d}{d\lambda})^m \frac{1}{1-b\lambda} = \frac{1}{1-b\lambda} \omega_m(\frac{b\lambda}{1-b\lambda}),$$

true for any $b$ independent of $\lambda$ (as follows from (6.2) by the substitution $\lambda \to b\lambda$).

**7. Evaluation of integrals.**

The family of polygamma functions $\psi^{(m)}$ is defined by: $\psi^{(0)} = \psi$ and

$$\psi^{(m)}(s) = \frac{d^m}{ds^m} \psi(s) = (-1)^{m+1} m! \zeta(m+1, s), \ m \in \mathbb{N}.$$

Polygamma functions of negative order, i.e for negative integers $m$, have been defined by several authors, as discussed in [1] and [7]. Adamchik [1], for example, considered the family of functions

$$\psi^{(-1)}(t) = \log\Gamma(t), \quad \psi^{(-k)}(t) = \frac{1}{(k-2)!} \int_0^t (t-s)^{k-2} \log\Gamma(s)\,ds, \quad k \geq 2, \tag{7.1}$$

and evaluated this integral in a closed form.

At the same time, it is clear from equation (4.3), i.e.

$$\sum_{n=1}^{\infty} \zeta(n+1, a) t^n = \psi(a) - \psi(a-t),$$

that evaluating the series (1.3) is equivalent to evaluating the integral:

$$\int_0^t s^{p-1} \psi(a-s)\,ds, \tag{7.2}$$

similar to the one in (7.1). For details, comments and other evaluations see [1], [7], [10], [11].

It is convenient to have at hand general method for evaluation of such integrals. The proofs of Theorems 1, 2, 3 show that when a function $f(a)$ has a representation by a Hankel integral as in (2.8), we can evaluate integrals of the form

$$\int_0^t s^m f(a+s)\,ds, \quad Re(a) > 0, \quad m = 0, 1, 2, \ldots$$

as combinations of the functions $g, \psi$ and polynomials. For example, we shall evaluate the integral (1.10). We need a slight modification of the integration rule (4.8). Namely:

$$\int_0^t y^m e^{zy}\,dy = e^{tz} \sum_{k=0}^{m} (-1)^k k! \binom{m}{k} \frac{t^{m-k}}{z^{k+1}} - (-1)^m \frac{m!}{z^{m+1}}, \tag{7.3}$$

for every $m = 0, 1, 2, \ldots$ .

**Theorem 4.** For $Re(a) > 0$, $|t| < Re(a)$ and $m = 0, 1, 2, \ldots$:

$$\int_0^t s^m \log\Gamma(a+s)\,ds = -\gamma \frac{t^{m+2}}{m+2} + [\log\sqrt{2\pi} - \gamma(a-\frac{1}{2})] \frac{t^{m+1}}{m+1} \tag{7.4}$$

$$+\sum_{k=0}^{m} \frac{(-1)^k}{k+1} \binom{m}{k} t^{m-k} g(k+1, a+t) - \frac{(-1)^m}{m+1} g(m+1, a).$$

The proof follows from Lemma 3, formula (7.3), and (2.8).

In particular, we have:

$$\int_0^t \log\Gamma(a+s)ds = -\gamma \frac{t^2}{2} + [\log\sqrt{2\pi} - \gamma(a-\frac{1}{2})]t + g(1, a+t) - g(1, a). \qquad (7.5)$$

Using (2.8) we write

$$g(1, a+t) - g(1, a) = \zeta'(-1, a+t) - \zeta'(-1, a) - \psi(2)[\zeta(-1, a+t) - \zeta(-1, a)].$$

Next, we shall simplify the expression

$$\psi(2)[\zeta(-1, a+t) - \zeta(-1, a)].$$

The second Bernoulli polynomial is

$$B_2(a) = a^2 - a + \frac{1}{6},$$

and according to (1.9):

$$\zeta(-1, a+t) - \zeta(-1, a) = \frac{-1}{2}(B_2(a+t) - B_2(a)) = \frac{-1}{2}(2at + t^2 - t).$$

Also, by (1.8)

$$\psi(2) = 1 - \gamma,$$

and we obtain from (7.5)

$$\int_0^t \log\Gamma(a+s)ds = -\frac{t^2}{2} + (\log\sqrt{2\pi} - a + \frac{1}{2})t + \zeta'(-1, a+t) - \zeta'(-1, a). \qquad (7.6)$$

The integral in (7.6) was first evaluated by Barnes [4] in terms of the $G$-function which he studied. The $G$-function, called now the Barnes function is defined by the difference equation

$$G(s+1) = \Gamma(s)G(s), \quad G(1) = 1.$$

In terms of $G(s)$ the evaluation is:

$$\int_0^t \log\Gamma(a+s)\,ds = -\frac{t^2}{2} + (\log\sqrt{2\pi} - a + \frac{1}{2})t \qquad (7.7)$$

$$+ (a+t-1)\log\Gamma(a+t) - \log G(a+t) - (a-1)\log\Gamma(a) + \log G(a),$$

([4], [19], p.207). The connection with (7.6) is easily found by using the equation ([19], p. 94):

$$\log G(a) = \frac{1}{12} - \log A + (a-1)\log\Gamma(a) - \zeta'(-1,a), \qquad (7.8)$$

where $A$ is the Glaisher-Kinkelin constant ([1], [6], [19]).

The function $\zeta'(-1,a)$ a has a Hankel integral representation. From (1.8) and (1.9):

$$\psi(2)\zeta(-1,a) = \frac{-1}{2}(1-\gamma)(a^2 - a + \frac{1}{6}),$$

and by (2.8):

$$\zeta'(-1,a) = \frac{1}{2}(1-\gamma)(a^2 - a + \frac{1}{6}) + \frac{1}{2\pi i}\int_L \frac{z^{-2}e^{az}}{1-e^z}\operatorname{Log} z\,dz. \qquad (7.9)$$

We can use equation (7.8) now to represent $\log G(a)$ by a Hankel integral, as $\zeta'(-1,a)$ and $\log\Gamma(a)$ already have such representations. Substituting (7.9) and (3.11) in (7.8) one arrives at:

$$\log G(a) = p(a) + \frac{1}{2\pi i}\int_L [(a-1)z^{-1} - z^{-2}]\frac{e^{az}\operatorname{Log} z}{1-e^z}\,dz, \qquad (7.10)$$

or

$$\log G(a) = p(a) + (a-1)g(0,a) - g(1,a), \qquad (7.11)$$

where $p(a)$ is the second degree polynomial:

$$p(a) = \frac{-1}{2}(1+\gamma)a^2 + (\log\sqrt{2\pi} + \gamma + \frac{1}{2})a - \frac{5\gamma}{12} - \log(A\sqrt{2\pi}). \qquad (7.12)$$

Therefore, one can evaluate integrals of the form

$$\int_0^t s^m \log G(a+s)\,ds, \quad Re(a) > 0, \quad m = 0, 1, 2, \ldots,$$

either by the method used in Theorems 1-4 or, alternatively, by repeated integration of (7.11), using the rule (2.10) for integrating $g(n,a)$.

**Appendix**

This is a list of some Hankel integral representations obtained in the paper. They are given with their original numbers. The contour $L$ is described at the beginning of Section 2.

$$\zeta(s,a) = \frac{\Gamma(1-s)}{2\pi i} \int_L \frac{z^{s-1}}{1-e^z} e^{az}\,dz; \tag{2.4}$$

$$\zeta(-n,a) = \frac{n!}{2\pi i} \int_L \frac{z^{-(n+1)}}{1-e^z} e^{az}\,dz, \quad n = 0, 1, \ldots; \tag{2.5}$$

$$\zeta(n+1,a) = \frac{(-1)^{n+1}}{2\pi i n!} \int_L \frac{z^n \operatorname{Log} z}{1-e^z} e^{az}\,dz, \quad n = 1, 2, \ldots; \tag{2.6}$$

$$g(n,a) \equiv \zeta'(-n,a) + \psi(n+1)\zeta(-n,a) = \frac{n!}{2\pi i}\int_L \frac{z^{-n-1} e^{az} \operatorname{Log} z}{1-e^z}\,dz \tag{2.8}$$

$(n = 0, 1, 2, \ldots);$

$$\psi(a) + \gamma = \frac{1}{2\pi i} \int_L \frac{e^{az}}{1-e^z} \operatorname{Log} z\,dz; \tag{3.7}$$

$$\log \Gamma(a) = \log\sqrt{2\pi} - \gamma\left(a - \frac{1}{2}\right) + \frac{1}{2\pi i} \int_L \frac{z^{-1} e^{az} \operatorname{Log} z}{1-e^z}\,dz; \tag{3.11}$$

$$\Phi(\lambda, s, a) = \frac{\Gamma(1-s)}{2\pi i} \int_L \frac{z^{s-1} e^{az}}{1-\lambda e^z}\,dz; \tag{5.2}$$

**References**.

1. **V. Adamchik,** PolyGamma functions of negative order, *J. Computational Applied Math.*, 100 (1998),191-199.

2. **G. E. Andrews, R. Askey, R. Roy,** *Special Functions*, Cambridge University Press, 1999.

3. **Tom. M. Apostol**, *Introduction to Analytic Number Theory*, Springer-Verlag, New York, 1986.

4. **E. W. Barnes**, The theory of the G-function, *Quart.J. Math.*, 31 (1899/1990), 264-314. (Originally *Quarterly Journal of Pure Mathematics*.)

5. **Khristo N. Boyadzhiev**, A series transformation formula and related polynomials, Int. J. Math. Math. Sci., 2005:23 (2005) 3849-3866.

6. **J. Choi, H. M. Srivastava,** Certain classes of series associated with the Zeta function and multiple Gamma functions, *J. Comput. Appl. Math*., 118 (2000), 87-109.

7. **J. Choi, H. M. Srivastava and V. S. Adamchik**, Multiple Gamma and related functions, *Appl. Math.Comput.*, 134 (2003), 515-533.

8. **A. Erdélyi** (editor), *Higher Transcendental Functions*, vol. 1, New York: McGrow-Hill, 1955

9. **O. Espinosa and V. Moll**, A generalized polygamma function, *Integral Transforms and Special Functions*, v. 15, no. 2, (2004), 101-115.

10. **O. Espinosa and V. Moll**, On some integrals involving the Hurwitz zeta function: Part 2, *The Ramanujan J.*, 6, no. 4 (2002) 449-468.

11. **R. Wm. Gosper, Jr.**, $\int_{n/4}^{m/6} \ln\Gamma(z)\,dz$, *Fields Ins. Communications*, 14 (1997), 71-76.

12. **M. Hashimoto, S. Kanemitsu, Y. Tanigawa, M. Yoshimoto, W.-P. Zhang**, On some slowly convergent series involving the Hurwitz zeta-function, *J. Computational Applied.*

k-boyadzhiev@onu.edu
Ohio Northern University
Department of Mathematics
Ada, OH 45810, USA